\documentclass{article}
\usepackage{amssymb}
\usepackage{amsfonts}
\usepackage{amsmath}

\begin{document}

\begin{center}
{\Large Special numbers, special quaternions and special symbol elements}

\begin{equation*}
\end{equation*}%
Diana SAVIN 
\begin{equation*}
\end{equation*}
\end{center}

\textbf{Abstract. }{\small In this paper we define and we study properties of} $\left( l,1,p+2q,q\cdot l\right) -$ {\small numbers,} $\left( l,1,p+2q,q\cdot l\right) -$ {\small quaternions,} $\left( l,1,p+2q,q\cdot l\right) -$ {\small symbol elements. Finally, we obtain 
an algebraic structure with these elements.} 

\medskip


\medskip

\begin{equation*}
\end{equation*}

\textbf{1. Introduction}%
\begin{equation*}
\end{equation*}%

Quaternion algebras and of symbol algebras have applications in various branches of mathematics, but also in computer science, physics, signal theory.

In this chapter we introduce special numbers, special quaternions, special symbol elements, and we present some of their properties and their applications in combinatorics, number theory and associative algebra theory.

Let $K$ be a field with $char\left(K\right)\neq2$ and let $\alpha, \beta$$\in$$K\backslash \{0\}.$ We recall that the generalized quaternion algebra $H_{K}\left(\alpha, \beta\right)$ is an algebra over the field $K$ with a basis $\{e_{1},e_{2},e_{3},e_{4}\}$ (where $e_{1}=1$ ) and the following multiplication:
\begin{equation*}
\begin{tabular}{c||c|c|c|c|}
$\cdot $ & $1$ & $e_{2}$ & $e_{3}$ & $e_{4}$ \\ \hline\hline
$1$ & $1$ & $e_{2}$ & $e_{3}$ & $e_{4}$ \\ \hline
$e_{2}$ & $e_{2}$ & $\alpha$ & $e_{4}$ & $\alpha e_{3}$ \\ \hline
$e_{3}$ & $e_{3}$ & $-e_{4}$ & $\beta$ & $-\beta e_{2}$ \\ \hline
$e_{4}$ & $e_{4}$ & $-\alpha e_{2}$ & $\beta e_{2}$ & $-\alpha
\beta$ \\ \hline
\end{tabular}%
\end{equation*}%
Let $x$ be element from $H_{K}\left(\alpha, \beta\right),$ $x=x_{1}\cdot 1+x_{2}e_{2}+x_{3}e_{3}+x_{4}e_{4},$
where $x_{i}\in K,$ $\left(\forall\right)$$i\in \{1,2,3,4\}$ and let $\overline{x}$ be the conjugate of $x,$ $x=x_{1}\cdot 1-x_{2}e_{2}-x_{3}e_{3}-x_{4}e_{4}.$
The trace of $x$ is $\boldsymbol{t}\left(
x\right) =x+\overline{x}=2x_{1}.$ The norm of $x$ is $\boldsymbol{n}\left(
x\right) =x\cdot \overline{x}=x_{1}^{2}-\alpha x_{2}^{2}-\beta x_{3}^{2}+\alpha \beta x_{4}^{2}.$%\\

If $K=\mathbb{R}$ and $\alpha=\beta=-1,$ we obtain Hamilton quaternion algebra $\mathbb{H}_{\mathbb{R}}\left( -1 ,-1 \right),$ with the basis $%
\{1,i,j,k\}.$

The generalization of a quaternion algebra is a symbol algebra. 

Let $n$ be an arbitrary positive integer, $n\geq 3$ and let $K$ be a field with $char(K),$ which does not divide $n,$ containing
$\xi,$ where $\xi$ is a primitive $n$-th root of unity. Let $a,b$ $\in K\backslash \{0\}.$ The algebra $A$ over $K$ generated by elements $x$ and $y$ where%
\begin{equation*}
x^{n}=a,y^{n}=b,yx=\xi xy
\end{equation*}
is called a \textit{symbol algebra} and it is denoted by $\left( \frac{a,~b}{%
K,\xi }\right).$ Symbol algebras are also known as \textit{%
power norm residue algebras}.
For $n = 2,$ we obtain the quaternion algebra over the field $K.$ Quaternion algebras and symbol algebras are associative but non-commutative algebras, of dimension $n^{2}$ over $K.$ Also, they are central simple over the field $K$ (this means they are simple algebras and their centers are equal to $K$). Theoretical aspects about these algebras can be found in the books [Pi; 82], [La;04],   [Mil], [Gi, Sz; 06], [Le; 05], [Als, Ba; 04], [Vi; 80], [Vo; 10], [Ko], [Mi; 71]. Several properties of these algebras and their applications in number theory, combinatorics, associative algebra, geometry, coding theory, mechanics can be found in the articles
[Ak, Ko, To; 14], [Fla; 12], [Fl, Sa, Io; 13], [Fl, Sa; 14], [Fl, Sa; 15], [Fl, Sa; 15 (1)], [Fl, Sa; 15 (2)], [Fl, Sa; 17], [Fl, Sa; 18], [Fl, Sh; 13], [Fl, Sh; 13 (1)], [Ha; 12], [Ho; 63], [Ja, Ya; 13], [Ka,Ha; 17], [Li; 12], [Ra; 15], [Sa, Fa, Ci; 09], [Sa; 14 (1)], [Sa; 16], [Sa; 16 (1)], [Sa; 17],  [Sa; 17 (1)], [Sw; 73], [Ta; 13]. 

Many mathematicians studied the Fibonacci numbers, Lucas numbers, Pell numbers, Pell-Lucas numbers, generalized Fibonacci-
Lucas numbers, the generalized Pell- Fibonacci-Lucas numbers, Fibonacci polynomials, Jacobsthal-Lucas polynomials, Fibonacci quaternions, the generalized Fibonacci-Lucas quaternions, the generalized Pell-Fibonacci-Lucas quaternions (see [Ho; 63], [Ca; 15],   [Ca, Mo; 16], [Ca; 16], [Fl, Sa; 15], [Fl, Sa; 18], [Fl, Sh; 13], [Ha; 12], [Sa; 14], [Sw; 73], etc.).

 In this paper we define $\left( l,1,p+2q,q\cdot l\right) -$ numbers, $\left( l,1,p+2q,q\cdot l\right) -$  quaternions, $\left( l,1,p+2q,q\cdot l\right) -$ symbol elements. We also study properties and applications of these elements.

This paper is organized as follow: section 2 is a preliminary section, containing theoretical notions
 which we will then use in our results. In Section 3 we introduce two special number sequences (namely $\left(a_{n}\right) _{n\geq 0},$
$\left(b_{n}\right) _{n\geq 0}$), we obtain some interesting properties of these sequences and we also obtain some quaternion algebras which split or some division quaternion algebras. After these, in the same section, we introduce $\left( l,1,p+2q,q\cdot l\right) -$ numbers, $\left( l,1,p+2q,q\cdot l\right) -$  quaternions, $\left( l,1,p+2q,q\cdot l\right) -$ symbol elements and we obtain interesting properties and applications of them.

\begin{equation*}
\end{equation*}

\textbf{2. Preliminaries }%
\begin{equation*}
\end{equation*}%
\qquad \qquad

First of all, we recall some results about prime integers, about diophantine equations or about the Fibonacci numbers, properties which will be necessary (in the next section) for to study some quaternion algebras.\newline
\smallskip\\
\textbf{Proposition 2.1.} ([Cu; 06]). \textit{Let} $m$ \textit{be a fixed positive integer.
The diophantine equation} $x^{2}+my^{2}=z^{2}$ \textit{has an infinity of solutions}:
\begin{equation*}
x=a^{2}- mb^{2}, y=2mab, z=a^{2}+ mb^{2}, a,b\in \mathbb{Z}.
\end{equation*}%
\textbf{Theorem 2.2.} ([Al, Go; 99]). \textit{Let n be a positive integer. Then, there exist integers} $x, y$ \textit{such that} $n = x^{2} + y^{2}$ \textit{if and only if the exponent of any prime} $p\equiv3$ (\textit{mod} $4$) \textit{that divides} $n$ \textit{is even}.\\
\smallskip\\
\textbf{Proposition 2.3.} ([Sa; 14]). \textit{For each positive integer} $n,$ $n\equiv7$ (\textit{mod} $16$), \textit{there exist integer numbers} $x, y$ \textit{so that, the Fibonacci number} $f_{n}$ \textit{can be written as} $f_{n} = x^{2} + 9y^{2}.$\\
\smallskip\\
Let $K$ be a field with $char\left(K\right)\neq2,$ let $\alpha, \beta$$\in$$K\backslash \{0\}$ and let  the generalized quaternion algebra $H_{K}\left(\alpha, \beta\right).$ $H_{K}\left(\alpha, \beta\right)$ is a division algebra if and only if for
$x\in$$H_{K}\left(\alpha, \beta\right)$ we have $\boldsymbol{n}\left(x\right)=0$ only for $x = 0.$\\
We recall that $H_{K}\left(\alpha, \beta\right)$ is called split by $K$ if it is isomorphic with a matrix algebra over $K$ (see [Pi; 82], [La;04], [Mil], [Gi, Sz; 06]). It is known the following remark about the quaternion algebras.\\
\smallskip\\
\textbf{Remark 2.4.} ([La;04], [Le; 05]). \textit{Let} $K$ \textit{be a field with} \textit{char}$K\neq 2$ \textit{and let} $\alpha, \beta$$\in$$K\backslash \{0\}.$ \textit{Then, the quaternion algebra} $H_{K}\left(\alpha, \beta\right)$ \textit{is either split or a division algebra.}\\
\smallskip\\
In the book [Gi, Sz; 06] appears the following criterion to decide if a quaternion algebra splits.\\
\smallskip\\
\textbf{Proposition 2.5.} ([Gi, Sz; 06]). \textit{Let} $K$ \textit{be a field with} \textit{char}$K\neq 2$ \textit{and let} $\alpha ,\beta$$\in$$K\backslash \{0\}.$ \textit{The quaternion algebra }$\mathbb{H}%
_{K}\left( \alpha ,\beta \right) $\textit{\ splits if and only if the
conic }$C\left( \alpha ,\beta \right) :$ $\alpha x^{2}+\beta y^{2}=z^{2}$ 
\textit{\ has a rational point over }$K($\textit{i.e. if there are }$%
x_{0},y_{0},z_{0}\in K,$\textit{not all zero such that }$\alpha x_{0}^{2}+\beta
y_{0}^{2}=z_{0}^{2}).$

\begin{equation*}
\end{equation*}%
\textbf{3. Some properties of special quaternions}%
\begin{equation*}
\end{equation*}

Let $l$ be a nonzero natural number. We consider the sequence $\left( a_{n}\right) _{n\geq 0}$ 
\begin{equation*}
a_{n}=l\cdot a_{n-1}+ a_{n-2},\;n\geq 2,a_{0}=0,a_{1}=1
\end{equation*}%
and let the sequence $\left( b_{n}\right) _{n\geq 0}$
\begin{equation*}
b_{n}=l\cdot b_{n-1}+b_{n-2},\;n\geq 2,b_{0}=2,b_{1}=l.
\end{equation*}%

Let $\alpha =\frac{l+\sqrt{l^{2}+4}}{2}$ and $\beta =\frac{l-\sqrt{l^{2}+4}}{2}.$ It results immediately the following relations: \smallskip 
\newline
\textbf{Binet's formula for the sequence} $\left( a_{n}\right) _{n\geq 0}.$ 
\begin{equation*}
a_{n}=\frac{\alpha ^{n}-\beta ^{n}}{\alpha -\beta }=\frac{\alpha ^{n}-\beta
^{n}}{\sqrt{l^{2}+4}},\ \ \left( \forall \right) n\in \mathbb{N}.
\end{equation*}%
\textbf{Binet's formula for the sequence} $\left( b_{n}\right) _{n\geq 0}.$ 
\begin{equation*}
b_{n}=\alpha ^{n}+\beta ^{n},\ \ \left( \forall \right) n\in \mathbb{N}.
\end{equation*}%
In the following, we show that the product of two elements belonging to the sequences $\left( a_{n}\right)
_{n\geq 0},$ $\left( b_{n}\right) _{n\geq 0}$ are transformed into sums of elements belonging to the same sequences. Also, we find another properties of these sequences.\newline
\smallskip\\
\textbf{Proposition 3.1.} \textit{Let} $\left( a_{n}\right) _{n\geq 0},$ $%
\left( b_{n}\right) _{n\geq 0}$ \textit{be the sequences previously defined.
Then, the following equalities are true}:\newline
i) 
\begin{equation*}
b_{n}b_{n+m}=b_{2n+m}+\left( -1\right) ^{n}b_{m},\ \left( \forall \right)
n,m\in \mathbb{N};
\end{equation*}%
ii) 
\begin{equation*}
a_{n}b_{n+m}=a_{2n+m}+\left( -1\right) ^{n+1}a_{m},\ \left( \forall \right)
n,m\in \mathbb{N};
\end{equation*}%
iii) 
\begin{equation*}
a_{n+m}b_{n}=a_{2n+m}+\left( -1\right) ^{n}a_{m},\ \left( \forall \right)
n,m\in \mathbb{N};
\end{equation*}%
iv) 
\begin{equation*}
a_{n}a_{n+m}=\frac{1}{l^{2}+4}%
\left[ b_{2n+m}+\left( -1\right) ^{n+1}b_{m}\right],\ \left( \forall \right) n,m\in \mathbb{N};
\end{equation*}%
v) 
\begin{equation*}
b_{n}+b_{n+2}=\left(l^{2}+4\right)\cdot a_{n+1},\ \left( \forall \right) n\in \mathbb{N};
\end{equation*}%
vi) 
\begin{equation*}
a^{2}_{n}+a^{2}_{n+1}=a_{2n+1},\ \left( \forall \right) n\in \mathbb{N};
\end{equation*}%
vii) 
\begin{equation*}
b^{2}_{n}+b^{2}_{n+1}=\left(l^{2}+4\right)\cdot a_{2n+1},\ \left( \forall \right) n\in \mathbb{N};
\end{equation*}%

\textbf{Proof.} Let $n,m$ be two positive integers. Applying Binet's formulae, we have:%
\newline
i) 
\begin{equation*}
b_{n}b_{n+m}=\left( \alpha ^{n}+\beta ^{n}\right) \cdot \left( \alpha ^{n+m}+\beta
^{n+m}\right)=
\end{equation*}%
\begin{equation*}
 =\alpha ^{2n+m}+\beta ^{2n+m}+\alpha ^{n}\beta ^{n}\left( \alpha ^{m}+\beta ^{m}\right)
=b_{2n+m}+\left( -1\right) ^{n}b_{m}.
\end{equation*}%
ii) 
\begin{equation*}
a_{n}b_{n+m}=\frac{\alpha ^{n}-\beta ^{n}}{\alpha -\beta }\left( \alpha
^{n+m}+\beta ^{n+m}\right) =
\end{equation*}%
\begin{equation*}
=\frac{\alpha ^{2n+m}-\beta ^{2n+m}}{\alpha -\beta }-\frac{\alpha ^{n}\beta
^{n}\left( \alpha ^{m}-\beta ^{m}\right) }{\alpha -\beta }=a_{2n+m}+\left(
-1\right) ^{n+1}a_{m}.
\end{equation*}%
iii) 
\begin{equation*}
a_{n+m}b_{n}=\frac{\alpha ^{n+m}-\beta ^{n+m}}{\alpha -\beta }\left( \alpha
^{n}+\beta ^{n}\right) =
\end{equation*}%
\begin{equation*}
=\frac{\alpha ^{2n+m}-\beta ^{2n+m}}{\alpha -\beta }+\frac{\alpha ^{n}\beta
^{n}\left( \alpha ^{m}-\beta ^{m}\right) }{\alpha -\beta }=a_{2n+m}+\left(
-1\right) ^{n}a_{m}.
\end{equation*}%
iv) 
\begin{equation*}
a_{n}a_{n+m}=\frac{\alpha ^{n}-\beta ^{n}}{\alpha -\beta }\cdot \frac{\alpha
^{n+m}-\beta ^{n+m}}{\alpha -\beta }=
\end{equation*}%
\begin{equation*}
=\frac{1}{\left(\alpha -\beta\right)^{2}}\cdot \left[\alpha ^{2n+m}+\beta ^{2n+m}-\alpha ^{n}\beta ^{n}\left( \alpha
^{m}+\beta ^{m}\right) \right]=
\end{equation*}%
\begin{equation*}
=\frac{1}{l^{2}+4}%
\left[ b_{2n+m}+\left( -1\right) ^{n+1}b_{m}\right] .
\end{equation*}%
v) \begin{equation*}
b_{n}+b_{n+2}=\alpha ^{n}+\beta ^{n}+\alpha ^{n+2}+\beta ^{n+2}=
\end{equation*}
\begin{equation*}
=\alpha ^{n+1}\cdot\left(\alpha +\frac{1}{\alpha}\right) + \beta ^{n+1}\cdot\left(\beta +\frac{1}{\beta}\right)=
\end{equation*}
\begin{equation*}
=\alpha ^{n+1}\cdot \sqrt{l^{2}+4} - \beta ^{n+1}\cdot \sqrt{l^{2}+4}=\left(l^{2}+4\right)\cdot a_{n+1}.
\end{equation*}
vi) Applying iv) for $m=0,$ we have:
\begin{equation*}
a^{2}_{n}+a^{2}_{n+1}=\frac{1}{l^{2}+4}%
\left[ b_{2n}+\left( -1\right) ^{n+1}b_{0}+b_{2n+2}+\left( -1\right) ^{n+2}b_{0}\right]=
\end{equation*}%
\begin{equation*}
=\frac{1}{l^{2}+4}\left[ b_{2n}+b_{2n+2}\right].
\end{equation*}%
Applying v) we obtain:
\begin{equation*}
a^{2}_{n}+a^{2}_{n+1}=a_{2n+1}.
\end{equation*}%
vii) Applying v) we have:
\begin{equation*}
b^{2}_{n}+b^{2}_{n+1}=\left( \alpha ^{n}+\beta ^{n}\right)^{2} + \left( \alpha ^{n+1}+\beta ^{n+1}\right)^{2}=
\end{equation*}%
\begin{equation*}
=\alpha ^{2n}+\beta ^{2n}+2\left(-1\right)^{n} + \alpha ^{2n+2}+\beta ^{2n+2}+2\left(-1\right)^{n+1}=
\end{equation*}%
\begin{equation*}
=b_{2n} + b_{2n+2}=\left(l^{2}+4\right)\cdot a_{2n+1}.
\end{equation*}%
$\Box \smallskip $\\
Let $(f_{n})_{n\geq 0}$ be the Fibonacci sequence and let $(l_{n})_{n\geq 0}$
be the Lucas sequence. There are well known the Cassini's identities for Fibonacci and Lucas numbers:
\begin{equation*}
f_{n+1}f_{n-1}-f_{n}^{2}=\left( -1\right) ^{n},\;\left( \forall \right)
\;n\in \mathbb{N}^{\ast },
\end{equation*}%
and 
\begin{equation*}
l_{n+1}l_{n-1}-l_{n}^{2}=5\cdot \left( -1\right) ^{n-1},\;\left( \forall
\right) \;n\in \mathbb{N}^{\ast },
\end{equation*}%
Now, we obtain similarly results for the the sequences $\left( a_{n}\right)
_{n\geq 0},$ $\left( b_{n}\right) _{n\geq 0}.$\\
\medskip\\
\textbf{Proposition 3.2.} \textit{Let} $\left( a_{n}\right) _{n\geq 0},$ $%
\left( b_{n}\right) _{n\geq 0}$ \textit{be the sequences previously defined.
Then, the following identities are true}:\newline
i) 
\begin{equation*}
a_{n+1}a_{n-1}-a_{n}^{2}=\left( -1\right) ^{n},\;\left( \forall \right)
\;n\in \mathbb{N}^{\ast };
\end{equation*}%
ii) 
\begin{equation*}
b_{n+1}b_{n-1}-b_{n}^{2}=\left( -1\right) ^{n}\cdot\left(l^{2}+4\right),\;\left( \forall \right)
\;n\in \mathbb{N}^{\ast }.
\end{equation*}%
\textbf{Proof.}
i) \begin{equation*}
a_{n+1}a_{n-1}-a_{n}^{2}=\frac{\alpha ^{n+1}-\beta^{n+1}}{\sqrt{l^{2}+4}}\cdot \frac{\alpha ^{n-1}-\beta^{n-1}}{\sqrt{l^{2}+4}}
-\frac{\left(\alpha ^{n}-\beta^{n}\right)^{2}}{l^{2}+4}=
\end{equation*}%
\begin{equation*}
=\frac{-\left(-1\right)^{n-1}\cdot\left(\alpha ^{2}+\beta^{2}\right)+2\left(-1\right)^{n}}{l^{2}+4}=\frac{\left(-1\right)^{n}\cdot\left(b_{2}+2\right)}{l^{2}+4}=\left(-1\right)^{n}.
\end{equation*}
ii) \begin{equation*}
b_{n+1}b_{n-1}-b_{n}^{2}=\left(\alpha ^{n+1}+\beta^{n+1}\right)\cdot\left(\alpha ^{n-1}+\beta^{n-1}\right)-\left(\alpha ^{n}+\beta^{n}\right)^{2}=
\end{equation*}
\begin{equation*}
=\left(-1\right)^{n-1}\cdot \left(\alpha ^{2}+\beta^{2}\right)-2\left(-1\right)^{n}=\left(-1\right)^{n-1}\cdot\left(b_{2}+2\right)=\left(-1\right)^{n-1}\cdot\left(l^{2}+4\right).
\end{equation*}
$\Box \smallskip $\\
\textbf{Proposition 3.3.} \textit{Let}  $%
\left( b_{n}\right) _{n\geq 0}$ \textit{be the sequence previously defined.
Then, the followings are true}:\newline
i) \textit{if} $l$ \textit{is even, then} $b_{n}$ \textit{is even} $\left(\forall\right)$ $n\in \mathbb{N};$\\
ii) \textit{if} $l$ \textit{is odd, then} $b_{n}$ \textit{is even if and only if} $n\equiv0$ (\textit{mod} $3$);\\
iii) \textit{if} $n\equiv0$ (\textit{mod} $6$), \textit{then} $b_{n-1}\cdot b_{n+1}$$\equiv$$3$ (\textit{mod} $4$);\\
iv) \textit{if} $n\equiv3$ (\textit{mod} $6$), \textit{then} $b_{n-1}\cdot b_{n+1}$$\equiv$$1$ (\textit{mod} $4$).\\
\smallskip\\
\textbf{Proof.} For i), ii), iii) and iv) the proof is immediate, using the principle of mathematics induction (after $n\in \mathbb{N}$).
$\Box \smallskip $\\
\textbf{Proposition 3.4.} \textit{Let} $\left( a_{n}\right) _{n\geq 0},$ $%
\left( b_{n}\right) _{n\geq 0}$ \textit{be the sequences previously defined.
Then, the followings are true}:\newline
i) \textit{The quaternion algebra}$\mathbb{H}_{\mathbb{Q%
}}\left(-1 , f_{2n+1}\right)$ \textit{splits}, $\left(\forall\right)$ $n\in \mathbb{N}^{\ast };$\\
ii) \textit{The quaternion algebra} $\mathbb{H}_{\mathbb{Q%
}}\left(-1 , 5f_{2n+1}\right)$ \textit{splits}, $\left(\forall\right)$ $n\in \mathbb{N}^{\ast };$\\
iii) \textit{The quaternion algebra} $\mathbb{H}_{\mathbb{Q%
}}\left(-1 , a_{2n+1}\right)$ \textit{splits}, $\left(\forall\right)$ $n\in \mathbb{N}^{\ast };$\\
iv) \textit{The quaternion algebra} $\mathbb{H}_{\mathbb{Q%
}}\left(-1 , \left(l^{2}+4\right)\cdot a_{2n+1}\right)$ \textit{splits}, $\left(\forall\right)$ $n\in \mathbb{N}^{\ast };$\\
v) \textit{The quaternion algebra} $\mathbb{H}_{\mathbb{Q%
}}\left(-1 , f_{2n+1}f_{2n-1}\right)$ \textit{splits}, $\left(\forall\right)$ $n\in \mathbb{N}^{\ast };$\\
vi) \textit{The quaternion algebra} $\mathbb{H}_{\mathbb{Q%
}}\left(-1 , a_{2n+1}a_{2n-1}\right)$ \textit{splits}, $\left(\forall\right)$ $n\in \mathbb{N}^{\ast };$\\
vii)\textit{The quaternion algebra} $\mathbb{H}_{\mathbb{Q%
}}\left(-1 , -b_{n+1}b_{n-1}\right)$ \textit{is a division algebra}, $\left(\forall\right)$ $n\in \mathbb{N}^{\ast };$\\
viii) \textit{The quaternion algebra} $\mathbb{H}_{\mathbb{Q%
}}\left(1 , b_{n+1}b_{n-1}\right)$ \textit{splits}, $\left(\forall\right)$ $n\in \mathbb{N}^{\ast };$\\
ix) \textit{If} $l$ \textit{is odd and} $n\equiv 0 $(\textit{mod} $6$), \textit{then the quaternion algebra} $\mathbb{H}_{\mathbb{Q%
}}\left(-1 , b_{n+1}b_{n-1}\right)$ \textit{is a division algebra}, $\left(\forall\right)$ $n\in \mathbb{N}^{\ast };$\\
x) \textit{If} $6\nmid n$ \textit{and the exponent of any prime} $p\equiv 3 $(\textit{mod} $4$) \textit{that divides} $b_{n+1}b_{n-1}$ \textit{is even, then the quaternion algebra} $\mathbb{H}_{\mathbb{Q%
}}\left(-1 , b_{n+1}b_{n-1}\right)$ \textit{splits}, $\left(\forall\right)$ $n\in \mathbb{N}^{\ast };$\\
xi) \textit{The quaternion algebra}$\mathbb{H}_{\mathbb{Q%
}}\left(-9 , f_{n}\right)$ \textit{splits}, $\left(\forall\right)$ $n\in \mathbb{N}^{\ast },$ $n\equiv 7$ (mod $16$).\\
\smallskip\\
\textbf{Proof.} Since iii) is a generalization of i), we are proving directly iii).\\
iii) If we consider the equation $-x^{2}+a_{2n+1}\cdot y^{2}=z^{2},$ we apply Proposition 3.1 (vi) and we obtain that it has the following solution
in $\mathbb{Q}\times \mathbb{Q}\times \mathbb{Q}:$ $\left(x_{0}, y_{0}, z_{0}\right)=\left(a_{n}, 1, a_{n+1}\right).$ According to Proposition 2.5, it results that the quaternion algebra $\mathbb{H}_{\mathbb{Q%
}}\left(-1 , a_{2n+1}\right)$ splits, $\left(\forall\right)$ $n\in \mathbb{N}^{\ast }.$\\
iv) Using  Proposition 3.1 (vii), it results that the equation $-x^{2}+\left(l^{2}+4\right)\cdot a_{2n+1}\cdot y^{2}=z^{2}$
has a solution in $\mathbb{Q}\times \mathbb{Q}\times \mathbb{Q},$ namely $\left(x_{0}, y_{0}, z_{0}\right)=\left(b_{n}, 1, b_{n+1}\right).$ Applying Proposition 2.5, it results that the quaternion algebra $\mathbb{H}_{\mathbb{Q%
}}\left(-1 , \left(l^{2}+4\right)\cdot a_{2n+1}\right)$ splits, $\left(\forall\right)$ $n\in \mathbb{N}^{\ast }.$\\
ii) This is a particular case of ii) (for $l=1$).\\
vi) This is a generalization of v), so we are proving only vi).\\
Using Proposition 3.2 i), we find the following solution in $\mathbb{Q}\times \mathbb{Q}\times \mathbb{Q}$ for the equation           $-x^{2}+a_{2n+1} \cdot a_{2n-1}\cdot y^{2}=z^{2}:$ $\left(x_{0}, y_{0}, z_{0}\right)=\left(\left(-1\right)^{n}, 1, a_{2n}\right).$ Applying Proposition 2.5, we obtain that the quaternion algebra $\mathbb{H}_{\mathbb{Q%
}}\left(-1 , a_{2n+1}a_{2n-1}\right)$ splits, $\left(\forall\right)$ $n\in \mathbb{N}^{\ast }.$\\
vii) Let the quaternion algebra $\mathbb{H}_{\mathbb{Q%
}}\left(-1 , -b_{n+1}b_{n-1}\right)$ and let $\{1,e_{2},e_{3},e_{4}\}$ a basis in this algebra.
Let $x=x_{1}\cdot 1 + x_{2}\cdot e_{2} + x_{3}\cdot e_{3}+ x_{4}\cdot e_{4}$$\in$$\mathbb{H}_{\mathbb{Q}}\left(-1 , -b_{2n+1}b_{2n-1}\right).$ The norm of $x$ is $n\left(x\right)=x\cdot \overline{x}=x^{2}_{1}+ x^{2}_{2} + 
b_{n-1}b_{n+1}x^{2}_{3}+ b_{n-1}b_{n+1}x^{2}_{4}.$ Since $b_{n}$$\in$$ \mathbb{N}^{\ast },$ for $\left(\forall\right)$ $n\in \mathbb{N}^{\ast },$ it results that $n\left(x\right)=0$ if and only if $x=0.$ So, $\mathbb{H}_{\mathbb{Q%
}}\left(-1 , -b_{n+1}b_{n-1}\right)$ is a division algebra for $\left(\forall\right)$ $n\in \mathbb{N}^{\ast }.$\\
Similarly, it results immediately that $\mathbb{H}_{\mathbb{Q%
}}\left(-1 , -f_{n+1}f_{n-1}\right),$ $\mathbb{H}_{\mathbb{Q%
}}\left(-1 , -a_{n+1}a_{n-1}\right),$ $\mathbb{H}_{\mathbb{Q%
}}\left(-1 , -l_{n+1}l_{n-1}\right)$ are division algebras  for $\left(\forall\right)$ $n\in \mathbb{N}^{\ast }.$\\
viii) We study if the equation $x^{2}+ b_{n+1}b_{n-1}\cdot y^{2}=z^{2}$ has rational solutions. Applying Proposition 2.1, it results that the equation $x^{2}+ b_{n+1}b_{n-1}\cdot y^{2}=z^{2}$ has solutions in integer numbers, so it has solutions in the set of rational numbers. Using Proposition 2.5, we obtain that the quaternion algebra $\mathbb{H}_{\mathbb{Q%
}}\left(1 , b_{n+1}b_{n-1}\right)$ splits.\\
ix) If $n\equiv 0 $(mod $6$), according to Proposition 3.3, $b_{n-1}\cdot b_{n+1}$$\equiv$$3$ (mod $4$). We study if the equation 
$-x^{2}+b_{n+1}b_{n-1}\cdot y^{2}=z^{2}$ has integer solutions. We suppose that this equation has a solution $\left(x_{0}, y_{0}, z_{0}\right)$$\in$$\mathbb{Z}\times \mathbb{Z}\times \mathbb{Z}\setminus\left\{\left(0, 0, 0\right)\right\},$ g.c.d$\left(x_{0}, y_{0}\right)=$
g.c.d$\left(y_{0}, z_{0}\right)=$g.c.d$\left(y_{0}, z_{0}\right)=1.$ We have: $b_{n+1}b_{n-1}\cdot y^{2}_{0}$$\equiv$$0$ or $3$ (mod $4$), but $x^{2}_{0}+ z^{2}_{0}$$\equiv$$1$  or $2$ (mod $4$), so we cannot have $b_{n+1}b_{n-1}\cdot y^{2}_{0}=x^{2}_{0}+ z^{2}_{0}.$
It results that the equation $-x^{2}+b_{n+1}b_{n-1}\cdot y^{2}=z^{2}$ does not have integer solutions. We obtain immediately that the equation $-x^{2}+b_{n+1}b_{n-1}\cdot y^{2}=z^{2}$ does not have  solutions in the set of rational numbers, so the quaternion algebra $\mathbb{H}_{\mathbb{Q%
}}\left(-1 , b_{n+1}b_{n-1}\right)$ does not split.\\
Applying Remark 2.4, we obtain that 
the quaternion algebra $\mathbb{H}_{\mathbb{Q%
}}\left(-1 , b_{n+1}b_{n-1}\right)$ is a division algebra, $\left(\forall\right)$ $n\in \mathbb{N}^{\ast }.$\\
x) Case $1:$ $l$ is odd.\\
If $n\equiv3$ (mod $6$), according to Proposition 3.3 (iv) $b_{n-1}\cdot b_{n+1}$$\equiv$$1$ (mod $4$). If $n\equiv1$ or $2$ or 
$4$ or $5$ (mod $6$), according to Proposition 3.3 (ii) $b_{n-1}\cdot b_{n+1}$ is even.\\
Case $2:$ $l$ is even, according to Proposition 3.3 (i) $b_{n-1}\cdot b_{n+1}$ is even.\\
In all these cases, it is possible to exist a prime $p\equiv 3 $(mod $4$) that divides $b_{n+1}b_{n-1}.$
If the exponent of any prime $p\equiv 3 $(mod $4$) that divides $b_{n+1}b_{n-1}$ is even, according to Theorem 2.2 there exist
integers $x_{0}; z_{0}$ such that $b_{n+1}b_{n-1} = x^{2}_{0} + z^{2}_{0}.$ This implies that $\left(x_{0}, 1, z_{0}\right)$ is a solution in integer numbers for the equation $-x^{2}+b_{n+1}b_{n-1}\cdot y^{2}=z^{2}$, so, according Proposition 2.5, the quaternion algebra $\mathbb{H}_{\mathbb{Q%
}}\left(-1 , b_{n+1}b_{n-1}\right)$ splits, $\left(\forall\right)$ $n\in \mathbb{N}^{\ast }.$\\
xi) Let $n$ be a positive integer number, $n\equiv 7$ (mod $16$). Using Proposition 2.3 we obtain that there are $x_{0},z_{0}$$\in$$\mathbb{Z}$ such that $\left(x_{0}, 1, z_{0}\right)$ is a solution of the equation $-9x^{2}+f_{n}\cdot y^{2}=z^{2}.$ Applying Proposition 2.5, we obtain that the quaternion algebra $\mathbb{H}_{\mathbb{Q%
}}\left(-9 , f_{n}\right)$ splits, $\left(\forall\right)$ $n\in \mathbb{N}^{\ast },$ $n\equiv 7$ (mod $16$).\\
$\Box \medskip $\\
Let $p,q$ be two arbitrary integers and $\left( a_{n}\right) _{n\geq 0},$ $\left( b_{n}\right) _{n\geq 0}$ are
the sequences previously defined. If $n$$\in$$N^{*},$ $a_{-n}=\left(-1\right)^{n+1}\cdot a_{n}.$\\
Let the sequence $\left(
u_{n}\right) _{n\geq 0},$ 
\begin{equation*}
u_{n+1}=pa_{n}+qb_{n+1},\;n\geq 0. 
\end{equation*}%
To avoid confusion, we will use the notation $u_{n}^{p,q}$ for $u_{n}.$\newline
We remark that $u_{n}=lu_{n-1}+u_{n-2},\ \left( \forall \right) n\in \mathbb{%
N},n\geq 2,$\newline
We calculate $u_{0}=pa_{-1}+qb_{0}=p+2q,$ $u_{1}=pa_{0}+qb_{1}=q\cdot l.$ We call the elements of the sequence $\left( u_{n}\right)
_{n\geq 0}$ the $\left( l,1,p+2q,q\cdot l\right) -$numbers.\\
\medskip\\
\textbf{Remark 3.5.} \textit{Let} $p,q$ \textit{be two arbitrary integers,
and let} $\left(u_{n}^{p,q}\right) _{n\geq 1}$ \textit{the sequence previously defined.
Then, we have:} 
\begin{equation*}
pa_{n+1}+qb_{n}=u_{n}^{p,q}+u_{n+1}^{pl,o},\forall ~n\in \mathbb{N}-\{0\}.
\end{equation*}%
\textbf{Proof.} We compute 
\begin{equation*}
pa_{n+1}+qb_{n}=pla_{n}+pa_{n-1}+qb_{n}=u_{n}^{p,q}+u_{n+1}^{pl,o}.
\end{equation*}%
$\Box \smallskip $\medskip \newline
Let $\alpha ,\beta $$\in $$\mathbb{Q}^{\ast }.$ We consider the generalized quaternion algebra $\mathbb{H}_{\mathbb{Q%
}}\left( \alpha ,\beta \right) $ with basis $\{1,e_{1},e_{2},e_{3}\}.$ We define the $n$%
-th $\left( l,1,p+2q,q\cdot l\right) -$ quaternion to be the element
of the form 
\begin{equation*}
U_{n}^{p,q}=u_{n}^{p,q}\cdot1+u_{n+1}^{p,q}\cdot e_{1}+u_{n+2}^{p,q}\cdot e_{2}+u_{n+3}^{p,q}\cdot e_{3}.
\end{equation*}%
\textbf{Remark 3.6. } \textit{Let} $U_{n}^{p,q}$ \textit{be the} $n$%
-th $\left( l,1,p+2q,q\cdot l\right) -$ \textit{quaternion. Then, we have}: 
\begin{equation*}
U_{n}^{p,q}=0\ \text{\textit{if\ and\ only\ if}}\ p=q=0.
\end{equation*}%
\textbf{Proof.} 
" $\Leftarrow $" It is trivial. \\
" $\Rightarrow $" If $U_{n}^{p,q}=0,$ using the fact that $%
\{1,e_{1},e_{2},e_{3}\}$ is a basis in quaternion algebra $%
\mathbb{H}_{\mathbb{Q}}\left( \alpha ,\beta \right) ,$ we obtain that $%
u_{n}^{p,q}=0, u_{n+1}^{p,q}=0, u_{n+2}^{p,q}=0, u_{n+3}^{p,q}=0.$\\
 From the recurrence relation of the sequence $\left( u_{n}^{p,q}\right) _{n\geq 1},$
it results that $u_{n-1}^{p,q}=0,$ $u_{n-2}^{p,q}=0,$ ..., $s_{1}^{p,q}=0$, $%
u_{0}^{p,q}=0.$ So, $q=0$ and $p=0.$\newline
$\Box \smallskip $\medskip \newline
About the  generalized
Fibonacci-Lucas quaternions $\left(G_{n}^{p,q}\right)_{n\geq0},$ in the paper [Fl, Sa; 15] (Theorem 3.5), we proved that:\\
i) \textit{The set}  
\begin{equation*}
M=\left\{ \sum\limits_{i=1}^{n}5G_{n_{i}}^{p_{i},q_{i}}|n\in \mathbb{N}^{\ast
},p_{i},q_{i}\in \mathbb{Z},(\forall )i=\overline{1,n}\right\} \cup \left\{
1\right\}
\end{equation*}%
\textit{has a ring structure with quaternion addition and multiplication}.\\
ii) \textit{The set} $M$ \textit{is an order of the quaternion algebra} $\mathbb{H}_{\mathbb{Q%
}}\left( \alpha ,\beta \right).$\\
We generalized these results for $\left(1, a, p+ 2q, q\right)$−quaternions $\left(S_{n}^{p,q}\right)_{n\geq0},$, in the paper [Fl, Sa; 17] (Proposition 5.4), namely:\\
\textit{Let} $a$ \textit{be a nonzero natural
number and let} $O$ \textit{be} \textit{the set} 
\begin{equation*}
O=\left\{ \sum\limits_{i=1}^{n}\left( 1+4a\right)
S_{n_{i}}^{p_{i},q_{i}}|n\in \mathbb{N}^{\ast },p_{i},q_{i}\in \mathbb{Z}%
,(\forall )i=\overline{1,n}\right\} \cup \left\{ 1\right\} .
\end{equation*}%
\textit{Then} $O$ \textit{is an order of the quaternion algebra} $\mathbb{H}%
_{\mathbb{Q}}\left( \alpha ,\beta \right).$\\
\smallskip\\
Similarly, in the paper [Fl, Sa; 18] we introduced the generalized Pell- Fibonacci-Lucas numbers $\left(r_{n}^{p,q}\right)_{n\geq0},$
 the generalized Pell- Fibonacci-Lucas quaternions $\left(R_{n}^{p,q}\right)_{n\geq0},$ and we proved
that (Proposition 3.7. from the paper [Fl, Sa; 18]) the set
\begin{equation*}
O=\left\{ \sum\limits_{i=1}^{n}8R_{n_{i}}^{p_{i},q_{i}}|n\in \mathbb{N}^{\ast
},p_{i},q_{i}\in \mathbb{Z},(\forall )i=\overline{1,n}\right\} \cup \left\{
1\right\}
\end{equation*}%
\textit{is an order of the quaternion algebra} $\mathbb{H}_{\mathbb{Q%
}}\left( \alpha ,\beta \right).$\\
\smallskip\\
Here, we generalized these numbers and these quaternions: the sequence $\left(a_{n}\right)_{n\geq0}$ is the generalization for 
the Pell sequence $\left(P_{n}\right)_{n\geq0}$ and the sequence $\left(b_{n}\right)_{n\geq0}$ is the generalization for 
the Pell-Lucas sequence $\left(Q_{n}\right)_{n\geq0}.$ Also, the sequence $\left(u_{n}^{p,q}\right)_{n\geq0}$ is the generalization for the sequence $\left(r_{n}^{p,q}\right)_{n\geq0}$ and the sequence of the $\left( l,1,p+2q,q\cdot l\right) -$ quaternions $\left(U_{n}^{p,q}\right)_{n\geq0}$ is the generalization for the sequence of the generalized Pell- Fibonacci-Lucas quaternions $\left(R_{n}^{p,q}\right)_{n\geq0}.$\\
\smallskip\\
Let $\epsilon$ be a primitive root of the unity of order $3$ and let $K$ be a field with the property $\epsilon$$\in$$K.$ Let 
$\alpha_{1}, \alpha_{2}$$\in$$K^{*}$ and let $A=\left(\frac{\alpha_{1}, \alpha_{2}}{K,\epsilon}\right)$ be the symbol algebra of degree $3.$ $A$ has a $K$- basis 
$\left\{x^{j_{1}}y^{j_{2}} |0\leq j_{1},j_{2}<3\right\},$ with $x^{3}=\alpha_{1},$ $y^{3}=\alpha_{2},$ $yx=\epsilon xy.$\\
In the paper [Fl, Sa; 14], we defined the $n$-th Fibonacci symbol element
$$F_{n} = f_{n}\cdot 1 + f_{n+1}\cdot x + f_{n+2}\cdot x^{2} + f_{n+3}\cdot y + f_{n+4}\cdot xy +$$
$$+ f_{n+5}\cdot x^{2}y + f_{n+6}\cdot y^{2}  + f_{n+7}\cdot xy^{2} + f_{n+8}\cdot x^{2}y^{2}.$$
In the paper [Fl, Sa, Io; 13] we defined the $n$-th Lucas symbol element
$$L_{n} = l_{n}\cdot 1 + l_{n+1}\cdot x + l_{n+2}\cdot x^{2} + l_{n+3}\cdot y + l_{n+4}\cdot xy
+$$
$$+ l_{n+5}\cdot x^{2}y + l_{n+6}\cdot y^{2}  + l_{n+7}\cdot xy^{2} + l_{n+8}\cdot x^{2}y^{2}.$$
Now, we define the $n$%
-th $\left( l,1,p+2q,q\cdot l\right) -$ symbol element to be the element
of the form 
\begin{equation*}
\mathbb{U}_{n}^{p,q}=u_{n}^{p,q}\cdot1+u_{n+1}^{p,q}\cdot x+u_{n+2}^{p,q}\cdot x^{2}+u_{n+3}^{p,q}\cdot y+
\end{equation*}%
\begin{equation*}
+u_{n+4}^{p,q}\cdot xy + u_{n+5}^{p,q}\cdot x^{2}y + u_{n+6}^{p,q}\cdot y^{2} + u_{n+7}^{p,q}\cdot xy^{2} + u_{n+8}^{p,q}\cdot x^{2}y^{2}.
\end{equation*}%
\textbf{Remark 3.7. } \textit{Let} $\mathbb{U}_{n}^{p,q}$ \textit{be the} $n$%
-th $\left( l,1,p+2q,q\cdot l\right) -$ \textit{symbol element. Then, we have}: 
\begin{equation*}
\mathbb{U}_{n}^{p,q}=0\ \text{\textit{if\ and\ only\ if}}\ p=q=0.
\end{equation*}%
The proof of this remark is similar to the proof of Remark 3.6.\\
\smallskip\\
With proof ideas similar to those in the Theorem 3.5 from the paper [Fl, Sa; 15], Proposition 5.4 from the paper [Fl, Sa; 17], Proposition 3.7. from the paper [Fl, Sa; 18], we obtain the following results:\\
\smallskip\\
\textbf{Proposition 3.8.} \textit{Let} $l$ \textit{be a nonzero natural
number and let} $M_{1}$ \textit{be} \textit{the set} 
\begin{equation*} 
M_{1}=\left\{ \sum\limits_{i=1}^{n}\left( l^{2}+4\right)
U_{n_{i}}^{p_{i},q_{i}}|n\in \mathbb{N}^{\ast },p_{i},q_{i}\in \mathbb{Z}%
,(\forall )i=\overline{1,n}\right\} \cup \left\{ 1\right\} .
\end{equation*}%
\textit{Then} $M_{1}$ \textit{is an order of the quaternion algebra} $\mathbb{H}%
_{\mathbb{Q}}\left( \alpha ,\beta \right) .$ \newline
\smallskip \newline
\textbf{Proposition 3.9.} \textit{Let} $l$ \textit{be a nonzero natural
number and let} $M_{2}$ \textit{be} \textit{the set} 
\begin{equation*} 
M_{2}=\left\{ \sum\limits_{i=1}^{n}\left( l^{2}+4\right)
\mathbb{U}_{n_{i}}^{p_{i},q_{i}}|n\in \mathbb{N}^{\ast },p_{i},q_{i}\in \mathbb{Z}%
,(\forall )i=\overline{1,n}\right\} \cup \left\{ 1\right\} .
\end{equation*}%
\textit{Then} $M_{2}$ \textit{is an order of the symbol algebra} $A=\left(\frac{\alpha_{1},\alpha_{2}}{K,\epsilon}\right).$ \newline
\smallskip \newline
Since the proofs of Proposition 3.8 and Proposition 3.9 are similar, we only prove one of them.\\
\smallskip \newline
\textbf{Proof of Proposition 3.9.}\\
We prove that $M_{2}$ is a free $\mathbb{Z}-$ submodule of rank $9$
of the symbol algebra $A=\left(\frac{\alpha_{1},\alpha_{2}}{K,\epsilon}\right).$\\
According to Remark 3.7, $\mathbb{U}_{n}^{0,0}=0$$\in $$O.$\newline
Let $n,m\in \mathbb{N}^{\ast },$ $p,q,p^{^{\prime }},q^{^{\prime }},c,d\in 
\mathbb{Z}.$ We have: 
\begin{equation*}
cu_{n}^{p,q}+du_{m}^{p^{^{\prime }},q^{^{\prime
}}}=u_{n}^{cp,cq}+u_{m}^{dp^{^{\prime }},dq^{^{\prime }}}.  
\end{equation*}%
This implies that
\begin{equation*}
c\mathbb{U}_{n}^{p,q}+d\mathbb{U}_{m}^{p^{^{\prime }},q^{^{\prime
}}}=\mathbb{U}_{n}^{cp,cq}+ \mathbb{U}_{m}^{dp^{^{\prime }},dq^{^{\prime }}} \tag{5.1.}
\end{equation*}%
So, $M_{2}$ is a free $\mathbb{Z}-$ submodule of rank $9$
of the symbol algebra $A.$\newline
We consider the set $M_{3}=\left\{ \sum\limits_{i=1}^{n}\left( l^{2}+4\right)
u_{n_{i}}^{p_{i},q_{i}}|n\in \mathbb{N}^{\ast },p_{i},q_{i}\in \mathbb{Z}%
,(\forall )i=\overline{1,n}\right\}.$
We are proving that $M_{2}$ is a subring of the symbol algebra $A.$ Taking into account the relation $\left(5.1.\right),$
 it is enough to prove that
$\left( l^{2}+4\right) \mathbb{U}_{n}^{p,q}\left( l^{2}+4\right) \mathbb{U}_{m}^{p^{^{\prime
}},q^{^{\prime }}}$$\in$$M_{2}.$ For this, it is enough to prove that 
$\left( l^{2}+4\right) u_{n}^{p,q}\left( l^{2}+4\right) u_{m}^{p^{^{\prime
}},q^{^{\prime }}}$$\in$$M_{3}.$\\
 Let $m,n$ be two integers, $n<m.$ We calculate:
\begin{equation*}
\left( l^{2}+4\right) u_{n}^{p,q}\left( l^{2}+4\right) u_{m}^{p^{^{\prime
}},q^{^{\prime }}}\text{=}
\end{equation*}%
\begin{equation*}
=\left( l^{2}+4\right) \left( pa_{n-1}+qb_{n}\right)
\left( l^{2}+4\right) \left( p^{^{\prime }}a_{m-1}+q^{^{\prime }}b_{m}\right) 
\text{=}
\end{equation*}%
\begin{equation*}
=\left(l^{2}+4\right) ^{2}pp^{^{\prime }}a_{n-1}a_{m-1}+\left(l^{2}+4\right)
^{2}pq^{^{\prime }}a_{n-1}b_{m}+
\end{equation*}%
\begin{equation*}
+\left(l^{2}+4\right) ^{2}p^{^{\prime }}qa_{m-1}b_{n}+\left(l^{2}+4\right)
^{2}qq^{^{\prime }}b_{n}b_{m}.
\end{equation*}%
Using Proposition 3.1, we have: 
\begin{equation*}
\left( l^{2}+4\right) u_{n}^{p,q}\cdot \left( l^{2}+4\right) u_{m}^{p^{^{\prime
}},q^{^{\prime }}}
\text{=}
\end{equation*}%
\begin{equation*}
=\left(l^{2}+4\right) ^{2}pp^{^{\prime }}\frac{1}{l^{2}+4%
}\left[ b_{n+m-2}+\left( -1\right) ^{n}b_{m-n}\right] +
\left(l^{2}+4\right) ^{2}pq^{^{\prime }}\left[ a_{m+n-1}+\left( -1\right)
^{n}a_{m-n+1}\right] +
\end{equation*}%
\begin{equation*}
\text{+}\left(l^{2}+4\right) ^{2}p^{^{\prime }}q\left[ a_{n+m-1}+\left(
-1\right) ^{n}a_{m-n-1}\right] \text{+}\left(l^{2}+4\right) ^{2}qq^{^{\prime }}%
\left[b_{n+m}+\left( -1\right) ^{n}b_{m-n}\right] \text{=}
\end{equation*}%
\begin{equation*}
\text{=}\left(l^{2}+4\right) ^{2}\left[ pq^{^{\prime }}a_{n+m-1}\text{+}%
qq^{^{\prime }}b_{m+n}\right] \text{+}\left(l^{2}+4\right) ^{2}\left[
\left( \text{-}1\right) ^{n}p^{^{\prime }}qa_{m-n-1}\text{+}\left( \text{-}%
1\right) ^{n}qq^{^{\prime }}b_{m-n}\right] \text{+}
\end{equation*}%
\begin{equation*}
+\left(l^{2}+4\right) \left[\left( -1\right) ^{n-1}\left(l^{2}+4\right)
pq^{^{\prime }}a_{m-n+1}+\left(-1\right) ^{n}pp^{^{\prime }}b_{m-n}\right] +
\end{equation*}%
\begin{equation*}
+\left(l^{2}+4\right) \left[ \left(l^{2}+4\right) p^{^{\prime
}}qa_{n+m-1}+pp^{^{\prime }}b_{n+m-2}\right]
\end{equation*}%
Applying the definition of the sequence $\left(u_{n}\right)_{n\geq0}$ and Remark 3.5, we obtain:
\begin{equation*}
\left( l^{2}+4\right) u_{n}^{p,q}\cdot \left( l^{2}+4\right) u_{m}^{p^{^{\prime
}},q^{^{\prime }}}\text{=}
\end{equation*}%
\begin{equation*}
\text{=}\left(l^{2}+4\right) u_{m+n}^{\left(l^{2}+4\right) pq^{^{\prime
}},\left(l^{2}+4\right) qq^{^{\prime }}}\text{+}\left(l^{2}+4\right)
u_{m-n}^{\left(-1\right) ^{n}\left(l^{2}+4\right) p^{^{\prime }}q,\left(
-1\right) ^{n}\left(l^{2}+4\right) qq^{^{\prime }}}\text{+}
\end{equation*}%
\begin{equation*}
+\left(l^{2}+4\right) u_{m-n}^{\left(-1\right) ^{n-1}\left(l^{2}+4\right)
pq^{^{\prime }},\left(-1\right) ^{n}pp^{^{\prime }}}+\left(l^{2}+4\right)
u_{m-n+1}^{\left(-1\right) ^{n-1}l\left(l^{2}+4\right) pq^{^{\prime }},0}+
\end{equation*}%
\begin{equation*}
+\left(l^{2}+4\right)u_{m+n-2}^{\left(l^{2}+4\right) p^{^{\prime
}}q,pp^{^{\prime }}}+\left(l^{2}+4\right) u_{m+n-1}^{l\left(l^{2}+4\right)
p^{^{\prime }}q,0}.
\end{equation*}%
So, $\left(l^{2}+4\right)u_{n}^{p,q}\cdot \left(l^{2}+4\right)
u_{m}^{p^{^{\prime }},q^{^{\prime }}}$$\in $$M_{3}.$\\
It results that $M_{2}$ is an order of the symbol algebra $A=\left(\frac{\alpha_{1},\alpha_{2}}{K,\epsilon}\right).$ 

$\Box \smallskip$
\begin{equation*}
\end{equation*}%
\textbf{References}\newline
\begin{equation*}
\end{equation*}%
[Ak, Ko, To; 14] M. Akyigit, HH. Koksal, M. Tosun, \textit{Fibonacci generalized quaternions}, Adv. Appl. Clifford Algebras 2014, 24(3), 631–-641.\newline
[Al, Go; 99] V. Alexandru, N.M. Gosoniu, \textit{Elements of Number Theory} (in Romanian), publisher Bucharest University, 1999.\newline
[Als, Ba; 04] M. Alsina, P. Bayer, \textit{Quaternion Orders, Quadratic Forms and Shimura Curves}, CRM Monograph Series,
22, American Mathematical Society, 2004. \newline
[Ca; 15] P. Catarino, \textit{A note on} $h\left(x\right)$ \textit{Fibonacci quaternion polynomials}, Chaos,
Solitons and Fractals, 77(2015), 1--5.\newline
[Ca, Mo; 16] P. Catarino, M. L. Morgado, \textit{On Generalized Jacobsthal and Jacobsthal-Lucas polynomials}, An. St. Univ. Ovidius Constanta, Mat. Ser., 24 (3) (2016), 61--78.\newline
[Ca; 16] P. Catarino,  \textit{The modified Pell and the modified k-Pell quaternions and octonions}, Adv. Appl. Clifford Algebras 26(2)(2016), 577–-590.\newline
[Cu; 06] I. Cucurezeanu, \textit{Equations in integer numbers} (in Romanian), publisher Aramis, 2006.\newline
[Fla; 12] R. Flatley, \textit{Trace forms of Symbol Algebras}, Alg. Colloquium, 19 (2012), 1117–-1124.\newline
[Fl, Sa, Io; 13] C. Flaut, D. Savin, G. Iorgulescu, \textit{Some properties of Fibonacci and Lucas symbol elements}, Journal of Mathematical Sciences: Advances and Applications, vol. 20, 2013, 37-43.\newline
[Fl, Sa; 14] C. Flaut, D. Savin, \textit{Some properties of symbol algebras of degree} $3$, Math. Reports
vol. 16 (66), no. 3 (2014), 443-463. \newline
[Fl, Sa; 15] C. Flaut, D. Savin, \textit{Quaternion Algebras and Generalized Fibonacci-Lucas Quaternions}, Adv. Appl. Clifford Algebras, 25(4)(2015),
853-862.\newline
[Fl, Sa; 15 (1)] C. Flaut, D. Savin, \textit{Some examples of division symbol algebras of degree 3 and 5}, Carpathian J. Math., 31 (2015), no. 2, 197–-204.\newline
[Fl, Sa; 15 (2)] C. Flaut, D. Savin,  \textit{About quaternion algebras and symbol algebras}, Bull. Univ.
Transilvania Brasov, Seria III, vol 7(56) (2014), no 2, 59--64.\newline
[Fl, Sa; 17] C. Flaut, D. Savin, \textit{Some remarks regarding} $a, b, x_{0}, x_{1}$ \textit{numbers and} $a, b, x_{0}, x_{1}$ \textit{quaternions,} submitted. \newline
[Fl, Sa; 18] C. Flaut, D. Savin, \textit{Some special number sequences obtained from a difference equation of degree three}, Chaos, Solitons and Fractals, 106 (2018), 67--71.\newline
[Fl, Sh; 13] C. Flaut, V. Shpakivskyi, \textit{On Generalized Fibonacci
Quaternions and Fibonacci-Narayana Quaternions}, Adv. Appl. Clifford
Algebras, 23(3)(2013), 673--688.\newline
[Fl, Sh; 13 (1)] C. Flaut, V. Shpakivskyi, \textit{Real matrix representations for the complex quaternions}, Adv. Appl. Clifford Algebras, 23(3) (2013), 657--671.\newline
[Gi, Sz; 06] Gille, P., Szamuely, T., \textit{Central Simple Algebras and Galois
Cohomology}, Cambridge University Press, 2006.\newline
[Ha; 12] S. Halici, \textit{On Fibonacci Quaternions}, Adv. in Appl.
Clifford Algebras, 22(2)(2012), 321--327.\newline
[Ho; 63] A. F. Horadam, \textit{Complex Fibonacci Numbers and Fibonacci
Quaternions}, Amer. Math. Monthly, 70(1963), 289--291.\newline
[Ja, Ya; 13] M.  Jafari, Y. Yayli, \textit{Rotation in four dimensions via Generalized Hamilton operators}, Kuwait J. Sci., 40 (2013),
No. 1, 67–-79.\newline
[Ka,Ha;17] A. Karatas, S. Halici, \textit{Horadam Octonions}, An St. Univ Ovidius Constanta, Mat Ser 2017; 25(3) (2017), 97--106.\newline
[Ko] D. Kohel, \textit{Quaternion algebras}, echidna.maths.usyd.edu.au/kohel/alg/doc/
AlgQuat.pdf\newline
[Ky; 15] I. Kyrchei, \textit{The column and row immanants over a split quaternion algebra}, Adv. Appl. Clifford Algebr, 25 (2015),
No. 3, 611–-619.\newline
[La; 04] T. Y. Lam, \textit{Introduction to Quadratic Forms over Fields}, AMS, 2004.\newline
[Le; 05] A. Ledet, \textit{Brauer Type Embedding Problems}, American Mathematical Society, 2005.\newline
[Li; 12] B. Linowitz, \textit{Selectivity in quaternion algebras}, Journal of Number Theory 132, 
1425-1437 (2012).\newline
[Mil] J.S. Milne, \textit{Class Field Theory}, http://www.math.lsa.umich.edu/ jmilne.\newline
[Mi; 71] J. Milnor, \textit{Introduction to Algebraic K-Theory}, Annals of Mathematics Studies, Princeton Univ. Press, 1971.\newline
[Pi; 82] R. S. Pierce, \textit{Associative Algebras}, Springer Verlag, 1982.\newline
[Ra; 15] J. L. Ramirez, \textit{Some Combinatorial Properties of the} $k$%
\textit{-Fibonacci and the} $k$\textit{-Lucas Quaternions}, An. St. Univ.
Ovidius Constanta, Mat. Ser., 23 (2)(2015), 201--212.\newline
[Sa, Fa, Ci; 09] D. Savin,  C. Flaut, C. Ciobanu, \textit{Some properties of the symbol algebras}, Carpathian Journal
of Mathematics , 25(2)(2009), 239--245\newline
[Sa; 14] D. Savin, \textit{Fibonacci primes of special forms}, Notes on Number Theory and Discrete
Mathematics, vol. 20, 2014, no.2, 10-19. \newline
[Sa; 14 (1)] D. Savin,  \textit{About some split central simple algebras}, An. St. Univ. Ovidius Constanta,
Mat. Ser., 22 (1) (2014), 263--272.\newline
[Sa; 16] D. Savin, \textit{About division quaternion algebras and division symbol algebras},
Carpathian Journal of Mathematics, vol. 32, No. 2 (2016), 233 -- 240.\newline
[Sa; 16 (1)] D. Savin, \textit{Quaternion algebras and symbol algebras over algebraic number field} $K,$
\textit{with the degree} $\left[K:\mathbb{Q}\right]$ \textit{even}, Gulf Journal of Mathematics, Vol 4, Issue 4 (2016), 16--
21.\newline
[Sa; 17] D. Savin, \textit{About Special Elements in Quaternion Algebras
Over Finite Fields}, Advances in Applied Clifford Algebras, June 2017, Vol.
27, Issue 2, 1801-1813.\newline
[Sa; 17 (1)] D. Savin, \textit{About split quaternion algebras over quadratic fields and symbol algebras
of degree n}, Bull. Math. Soc. Sci. Math. Roumanie, Tome 60 (108) No. 3, 2017, 307- 312.\newline
[Sw; 73] M. N. S. Swamy, \textit{On generalized Fibonacci Quaternions,} The
Fibonacci Quaterly, 11(5)(1973), 547--549.\newline
[Ta; 13] M. Tarnauceanu, \textit{A characterization of the quaternion group}, An. St. Univ. Ovidius Constanta, 21 (2013),
No. 1, 209–-214. \newline
[Vi; 80] Vigneras, M.F., \textit{Arithmetique des alg‘ebres de quaternions}, Lecture Notes in Math., no.
800, Springer, 1980.\newline
[Vo; 10] Voight, J., The Arithmetic of Quaternion Algebras, available on the author’s website:
http://www.math.dartmouth.edu/ jvoight/ crmquat/book/quat-modforms-041310.pdf,
2010.\newline

\begin{equation*}
\end{equation*}

\bigskip

Diana SAVIN

{\small Faculty of Mathematics and Computer Science, }

{\small Ovidius University, }

{\small Bd. Mamaia 124, 900527, CONSTANTA, ROMANIA }

{\small http://www.univ-ovidius.ro/math/}

{\small e-mail: \ savin.diana@univ-ovidius.ro, \ dianet72@yahoo.com}\bigskip
\bigskip

\end{document}